\documentclass[12pt,a4paper,reqno]{amsart}

\usepackage{anysize}
\usepackage[utf8]{inputenc}
\usepackage{amsmath}
\usepackage{amssymb}
\usepackage{array}
\usepackage[unicode]{hyperref}

\marginsize{2.5cm}{2.5cm}{3cm}{3cm}

\newtheorem{theorem}{Theorem}[section]

\begin{document}

\title{The WZ method and flawless WZ pairs}

\author{Jesús Guillera}
\address{Department of Mathematics, University of Zaragoza, 50009 Zaragoza, SPAIN}
\email{}
\date{}

\dedicatory{Dedicated to Doron Zeilberger on his $75^{th}$ birthday}

\begin{abstract}
Recently, Kam Cheong Au discovered a powerful methodology of finding new Wilf-Zeilberger (WZ) pairs. 
He calls it WZ seeds and gives numerous examples of applications to proving longstanding conjectural identities for reciprocal powers of $\pi$ and their duals for Dirichlet $L$-values. 
In this note we explain how a modification of Au's WZ pairs together with a classical analytic argument allows one to obtain simpler proofs of his results.
We illustrate our method with a few examples elaborated with assistance of Maple code that we have developed. 
\end{abstract}

\maketitle

\section{Wilf-Zeilberger (WZ) pairs}
Herbert Wilf and Doron Zeilberger invented the concept of WZ pair: Two hypergeometric (in $n$ and $k$) terms $F(n,k)$ and $G(n,k)$ form a WZ pair if the identity
\[
F(n+1,k)-F(n,k)=G(n,k+1)-G(n,k)
\]
holds.
A Maple code written by Zeilberger, available in Maple in the package 
\[
\texttt{with(SumTools[Hypergeometric]);}
\]
finds the mate of a term (that forms a WZ pair with it), whenever such exists, by means of a rational certificate $C(n,k)$ so that $G(n,k)=C(n,k) F(n,k)$ \cite{certificate}.

The simple concept of WZ pair has very useful and important properties. We display two of them.

\begin{theorem}
\label{th1.1}
Assume that $F(0,k)=0$ and $F(+\infty,k)=\lim_{n\to+\infty}F(n,k)=0$ for any complex $k$. Then
\[
\sum_{n=0}^{\infty} G(n,k+1)=\sum_{n=0}^{\infty} G(n,k) \quad\text{for all}\; k \in \mathbb{C}.
\]
In particular, the function
\[
g(k)=\sum_{n=0}^{\infty} G(n,k), 
\]
is periodic of period $k=1$.
\end{theorem}

Due to its periodicity the function $g(k)$ cannot have poles: those will be also (periodic) poles of the hypergeometric terms $F(n,k)$ which is impossible. However, there could be removable singularities.

\begin{proof}
For $k$ a complex number, we have
\begin{align*}
\sum_{n=0}^{\infty} (F(n+1,k)-F(n,k)) &=F(+\infty,k)-F(0,k) =0 \\ 
&=\sum_{n=0}^{\infty} (G(n,k+1)-G(n,k)) =g(k+1)-g(k).
\qedhere
\end{align*}
\end{proof}

We call WZ pairs meeting the conditions $F(0,k)=F(+\infty,k)=0$ for all $k\in\mathbb C$ flawless. It occurs very often that the hypotheses of Carlson's theorem hold and then we obtain a stronger result:
\[
g(k)=\sum_{n=0}^{\infty} G(n,k)=\text{constant} \quad\text{for all}\; k \in \mathbb{C}.
\]
For justification of the latter we use the following `periodic version' of \textbf{Carlson's theorem} \cite[Appendix]{Almkvist}: \\
\emph{Let $k=x+iy$, where $x$ and $y$ are the real and imaginary parts of $k$, respectively. If $g(k)$ is an holomorphic function such that $g(k)=g(k+1)$ and with some real number $c<2\pi$ we have $g(k)=\mathrm{O}(e^{c|y|})$ for $|y|$ sufficiently large, then $g(k)$ is a constant function}.

\begin{theorem}
For every WZ pair we have
\[
\sum_{n=0}^{\infty} G(n,0)-\lim_{k \to \infty} \sum_{n=0}^k G(n,k)
=\sum_{n=0}^{\infty} F(0,k)-\lim_{n \to \infty} \sum_{k=0}^n F(n,k),
\]
whenever both sides converge.
\end{theorem}

The theorem is used by Au in \cite{Pisco-seeds}.

\section{Transformations}
From a WZ pair we can get others; this mechanism forms the grounds of the WZ-seed methodology of Au explored in \cite{Pisco-seeds,Pisco-seeds-q}. The following transformations will be useful for our needs:
\begin{enumerate}
\item Transformations $F(n,k) \to F(sn, k \pm tn)$ with $s,t\in\mathbb Z_{>0}$ change the hypergeometric series and could produce acceleration or deceleration of its convergence rate.
\item Transformations $G(n,k) \to G(n \pm tk, sk)$ with $s,t\in\mathbb Z_{>0}$ obviously do not alter the series at $k=0$ but some of the new WZ pairs may have better properties than the original one. 
\item The \emph{duality} corresponds to the transformation $G(n,k) \to G(-n,-k)$ \cite[Section 4]{rama-divergent}. Observe that it converts divergent series into convergent ones and vice versa.
\end{enumerate}
Usually WZ pairs are not flawless (see \cite{A-equal-B} and \cite{Pisco-seeds}), but can be sometimes converted to such via suitable transformations.

\section{Ramanujan-like series and their duals}
In examples of this paper, we only consider rational Ramanujan-like series or their duals.

\subsection{Rational Ramanujan-like series}
These can be written in the form
\[
\sum_{n=0}^{\infty} R(n)
=\sum_{n=0}^{\infty} \left( \prod_{i=0}^{2m} \frac{(s_i)_n}{(1)_n} \right) z^n \sum_{k=0}^m a_k n^k 
= \frac{\sqrt{(-1)^m \chi}}{\pi^m},
\]
where $(a)_0=1$ and $(a)_n=a(a+1) \cdots (a+n-1)$ for $n \in \mathbb{Z}_{>0}$, is the rising factorial or Pochhammer symbol, $z$ is rational, $a_0, a_1,\dots,a_m$ are positive rationals, and $\chi$ is the discriminant of a certain quadratic field (imaginary or real), which is an integer.  
Let $z=z_0$, if $|z_0|>1$ then the series is divergent, but we understand it as analytic continuation of the series with $z$ at $z=z_0$, and we can evaluate them with Maple executing
\[
\texttt{sum(R(n),n=0..infinity,formal);}
\]
which converts them in hypergeometric form.

The series of degree $m=1$ were discovered and proved using modular equations using a method highlighted by the genius Srinivasa Ramanujan. In 2002--2003 we discovered several series of higher degree \cite{about-new-rama}. In addition, inspired by \cite{first-rama-wz}, we managed to prove with the WZ method the following famous Ramanujan formula for $1/\pi$ \cite{generators-rama}:
\[
\sum_{n=0}^{\infty} \frac{(\frac12)_n^3}{(1)_n^3} (42n+5) 
\left( \frac{1}{64} \right)^n=\frac{16}{\pi}.
\]
In addition, we discovered some new Ramanujan-like series of degree $2$ \cite{about-new-rama}, and proved the following one \cite{generators-rama}:
\[
\sum_{n=0}^{\infty} \frac{(\frac12)_n^5}{(1)_n^5} (820n^2+180n+13) 
\left( - \frac{1}{1024} \right)^n=\frac{128}{\pi^2},
\]
with the same strategy of the WZ method.

\subsection{Their duals}
The dual of a Ramanujan-like series for a reciprocal power of $\pi$ is obtained by replacing $n$ with $-n$ in accordance with the following dictionary:
\[
(1)_{-n} \to \frac{n (-1)^n}{(1)_n} \quad \text{and} \quad (a)_{-n} \to \frac{(-1)^n}{(1-a)_n}
\;\text{for}\; 0<a<1.
\]
If the sum of a rational Ramanujan-like series equals $\sqrt{(-1)^m \chi}/\pi^m$, then the sum of its dual is a rational multiple of the Dirichlet value $L_{\chi}(m+1)$ \cite{heuristic-supercong}. Some examples of them are $L_{-4}(2)=L((\frac{-4}{\cdot}),2)$ (Catalan's constant), $L_1(3)=\zeta(3)$, $L_1(5)=\zeta(5)$ and $L_1(4)=\zeta(4)$ (a rational multiple of $\pi^4$).
\par To summarize, the dual of a Ramanujan-like formula has the following form:
\[
\sum_{n=1}^{\infty} \left( \prod_{i=0}^{2m} \frac{(1)_n}{(s_i)_n} \right) \frac{z^{-n}}{(-n)^{m+1}} \sum_{k=0}^m a_k (-n)^{k}=r L_{\chi}(m+1),
\]
The series converges at $z=z_0$ if $|z_0|>1$, and analytic continuation is required if $|z_0|<1$, and we can determine the rational $r$ using supercongruences \cite{heuristic-supercong}. 
\par Finding flawless WZ pairs for these duals seems to be a challenging task.

\section{Examples inspired by Au's pairs}

\subsection{Proof of a Ramanujan-like formula of degree 4 conjectured by J.~Cullen}
We write the WZ pair of \cite[Example I]{Pisco-seeds} in a different way. Namely, we start from the entry
\[
U(n,k)=\frac{(\frac12)_n^5 (\frac12+k)_n^5 (\frac12 -k)_n (1+\frac{k}{4})_n (\frac14+\frac{k}{4})_n (\frac12+\frac{k}{4})_n (\frac34+\frac{k}{4})_n}{(1)_n^5 (1+\frac{k}{2})_n^5 (\frac12 + \frac{k}{2})_n^5} \frac{(\frac12)_k^4}{(1)_k^4},
\]
which involves only Pochhammer symbols of the forms $(a+bk)_n$ and $(c)_k$, where $a$ and $c$ are rationals in the interval $(0,1]$, take
\[
z=\left(\frac12\right)^{12}, \quad y=1,
\]
and construct the rational functions
\[
S(n,k)=2^{16} \frac{n^5 (4k + 6n + 1)(3 + 2k)(1 + 2k)}{(4n + k)(4k^2 + 8 k + 3)(2k + 1 - 2n)}
\]
and $R(n,k)$ using the recipe of \cite{Maple-on-WZ-pairs}.
Then the function $F(n,k)=U(n,k) S(n,k) z^n y^k$ has the WZ mate $G(n,k)=U(n,k) R(n,k) z^n y^k$. The conditions $F(0,k)=F(+\infty,k)=0$ of Theorem~\ref{th1.1} are satisfied, and we deduce that
\[
\sum_{n=0}^{\infty} G(n,k)=\sum_{n=0}^{\infty} G(n,k+1) \quad \text{for all}\; k \in \mathbb{C}.
\]
With the help of Carlson's theorem we arrive at the following stronger result:
\[
\sum_{n=0}^{\infty} G(n,k) = \text{constant} \quad \text{for all}\; k \in \mathbb{C}.
\]
To determine the value of the constant we take $k=1/2$ and notice that then we obtain the factor $(0)_n$ in the numerator of $U(n,1/2)$, hence in the numerator of $G(n,1/2)$. As $(0)_n=0$ if $n$ is different from $0$, the only nonzero term in the sum $\sum_{n=0}^{\infty} G(n,1/2)$ is $G(0,1/2)$, and we obtain
\[
\sum_{n=0}^{\infty} G(n,k)=\sum_{n=0}^{\infty} G\left(n,\frac12\right)=G \left(0,\frac12 \right)=\frac{2048}{\pi^4}.
\]
Thus, for $k=0$ we get
\[
\sum_{n=0}^{\infty} \frac{(\frac12)_n^7 (\frac14)_n (\frac34)_n}{n!^9} \,\frac{1}{2^{12 n}}(43680n^4 +20632n^3+4340n^2+466n+21)=\frac{2048}{\pi^4}.
\]
The Maple code in \cite{Maple-on-WZ-pairs} generalizes the pattern conjectured in \cite{on-WZ-pairs} for degree $1$ to Ramanujan-like series of higher degrees.
Finally, we notice that replacing $F(n,k)$ with $F(n,k-n)$ in the above leads to \cite[Example II]{Pisco-seeds}.

\subsection{Proof of one of our Ramanujan-like conjectured formulas of degree 2}
We write the WZ pair of \cite[Example VII]{Pisco-seeds} starting with
\[
U(n,k)=\frac{(\frac12)_n^3 (1+\frac{k}{3})_n (\frac13+\frac{k}{3})_n (\frac23+\frac{k}{3})_n (\frac16+\frac{2k}{3})_n (\frac12+\frac{2k}{3})_n (\frac56+\frac{2k}{3})_n }{(1)_n^3 (1+\frac{k}{2})_n^3 (\frac12+\frac{k}{2})_n^3} \frac{(\frac14)_k (\frac34)_k}{(1)_k^2},
\]
the values
\[
z=-\frac{3^6}{4^6}, \quad y=1,
\]
and constructing the rational functions
\[
S(n,k)=\frac{768 n^3}{k + 3n}
\]
and $R(n,k)$ as explained in \cite{Maple-on-WZ-pairs}.
Then the function $F(n,k)=U(n,k) S(n,k) z^n y^k$ has the WZ mate $G(n,k)=U(n,k) R(n,k) z^n y^k$. As in our first example, we deduce that
\[
\sum_{n=0}^{\infty} G(n,k) = \text{constant} \quad \text{for all}\; k \in \mathbb{C}. 
\]
To determine the constant we let $k \to +\infty$, taking into account that $(a+bk)_n \sim (bk)^n$, and we obtain the identity that specializes to
\[
\sum_{n=0}^{\infty} \frac{  (\frac12)_n (\frac13)_n (\frac23)_n (\frac16)_n (\frac56)_n }{n!^5} \left(-\frac{3^6}{4^6}\right)^n (1930n^2+549n+45) = \frac{384}{\pi^2}.
\]

\subsection{Proof of a formula for $\zeta(4)$}

Inspired by \cite[Example IV]{Pisco-seeds}, we write
\[
U(n,k)=\frac{(1)_n^3 (1+k)_n^3(1+2k)_n\left( 1+\frac{k}{2}\right)_n\left(\frac12+\frac{k}{2}\right)_n} {\left( \frac12
\right)_n^3 \left( \frac12+k \right)_n^3 \left(1+ \frac{2k}{3} \right)_n \left( \frac13+\frac{2k}{3} \right)_n \left( \frac23+\frac{2k}{3} \right)_n} \frac{(1)_k^4}{\left( \frac12\right)_k^4}
\]
and
\[
F(n,k)=U(n,k)S(n,k) z^n y^k,
\]
where $z=16/27$, $y=1$ and
\[
S(n,k)=\frac{n\left( 2n+2k+\frac32\right)}{(1+2k)(2n+2k+1)^3(3n+2k+1)(3n+2k+2)}.
\]
The function $F(n,k)$ has a WZ companion $G(n,k)$ such that 
\[
G(n,k)=U(n,k)R(n,k)z^n y^k,
\]
where $R(n,k)$ is a certain rational function that a Maple code written by Zeilberger determines.

On the other hand, as $F(n,k)$ vanishes for $n=0$ and as $n \to +\infty$, we deduce that
\[
\sum_{n=0}^{\infty} G(n,k) = \frac{\pi^4}{32} \quad \text{for all}\; k \in \mathbb{C},
\]
where we have determined the constant value of the sum by letting $k \to +\infty$. This includes the simpler particular case corresponding to $k=0$ \cite[Example IV]{Pisco-seeds}.

\subsection{Proof of a formula for $\zeta(2)$}

Playing with the term $F(n,k)$ of the WZ pair used in Example 3, and modifying some of the powers we discover a new WZ pair with
\[
U(n,k)=\frac{(1)_n (1+k)_n (1+2k)_n\left( 1+\frac{k}{2}\right)_n\left(\frac12+\frac{k}{2}\right)_n} {\left( \frac12 \right)_n \left( \frac12+k \right)_n \left(1+ \frac{2k}{3} \right)_n \left( \frac13+\frac{2k}{3} \right)_n \left( \frac23+\frac{2k}{3} \right)_n} \frac{(1)_k^2}{\left( \frac12\right)_k^2}
\]
and
\[
F(n,k)=U(n,k)S(n,k) z^n y^k,
\]
where $z=4/27$, $y=1$ and
\[
S(n,k)=\frac{n\left( 2n+2k+\frac32\right)}{(1+2k)(2n+2k+1)(3n+2k+1)(3n+2k+2)}.
\]
The function $F(n,k)$ has a WZ companion $G(n,k)$ such that 
\[
G(n,k)=U(n,k)R(n,k)z^n y^k,
\]
where $R(n,k)$ is the following rational function:
\[
R(n,k)=\frac{22n^3 + (34k + 49)n^2 + (12k^2 + 48k + 35)n + 8(k + 1)^2}{2(1 + 2n)(3 + 2k + 3n)(2n + 2k + 1)(3n + 2k + 1)(3n + 2k + 2)},
\]
which is discovered automatically by the Maple code.
On the other hand, as $F(n,k)$ vanishes for $n=0$ and as $n \to +\infty$, we deduce that
\[
\sum_{n=0}^{\infty} G(n,k) = \frac{\pi^2}{8} \quad \text{for all}\; k \in \mathbb{C},
\]
where we have determined the constant value of the sum by letting $k \to +\infty$. This includes the simpler particular case corresponding to $k=0$ \cite[eq. 14]{rama-divergent}

\subsection{Proof of the famous Gourevitch conjectured formula}
In the Appendix B of \cite{Pisco-seeds-q}, Au finds a WZ pair suitable for proving the Ramanujan-like series for $1/\pi^3$ conjectured by Gourevitch. His proof is somewhat complicated. The main obstacle there is that 
\[
F(0,k) \ne 0 \quad\text{and}\quad \lim_{n \to \infty} \sum_{k=0}^{\infty} F(n,k) \ne 0.
\]
The issue can be resolved by passing to $\tilde G(n,k)=G(n-k,k)$ in place of $G(n,k)$.
This leads to the following beautiful Ramanujan-like generalization:
\begin{multline*}
\sum_{n=0}^{\infty} \frac{(\frac12+3k)_n^2 (\frac12-k)_n^2 (\frac12+k)_n^3}{(1)_n^2(1+2k)_n^2(1+k)_n^3} \frac{(\frac12)_k (\frac16)_k^2 (\frac56)_k^2}{(1)_k^5} \left(\frac{1}{64}\right)^n \left(\frac{729}{1024}\right)^k
\\ \times 
(6n+6k+1)(28n^2+8n+1+12k^2+56nk+8k) = \frac{32}{\pi^3},
\end{multline*}
valid for all complex values of $k$ thanks to Carlson's theorem
(the value of the constant is obtained by taking $k=1/2$).
This (WZ provable) generalization proves the Gourevitch formula by specicializing $k=0$, namely,
\[
\sum_{n=0}^{\infty} \frac{(\frac12)_n^7}{(1)_n^7} \left(\frac{1}{64}\right)^n 
(6n+1)(28n^2+8n+1) = \frac{32}{\pi^3}.
\]
Full details of this derivation can be found in \cite{Maple-on-WZ-pairs}.

\section{Formulas with $|z|>1$}

If $|z|>1$ then hypergeometric series in the formulas diverge. In 1908, Barnes showed that the complex integral representation
\[
\sum_{n=0}^{\infty} \frac{(a)_n(b)_n}{(1)_n(c)_n} z^n
= \frac{1}{2 \pi i} \int_{-i \infty}^{i \infty} \frac{(a)_s (b)_s}{(c)_s} \Gamma(-s) ({-z})^s \,ds
\]
valid for the series within the domain $|z|<1$, $z\notin[0,-1]$, gives its analytic continuation to $z\in\mathbb C\setminus[0,\infty)$.
Here the integration path connecting $-i\infty$ with $i\infty$ is chosen in such a way that the poles of $\Gamma(-s)$ are all to the left from the contour, while the poles of $(a)_s (b)_s=\Gamma(a+s)\Gamma(b+s)/(\Gamma(a)\Gamma(b))$ are all to the right from it; the branch of $(-z)^s=\exp(s\ln(-z))$ is chosen to be real valued for $z<0$.

Inspired by the WZ pair in \cite[Example XII]{Pisco-seeds} we take its dual and write the new function $F(n,k)$ as $F(n,k)=U(n,k) S(n,k) z^n$, where $z=-3125/1024$, 
\[
U(n,k)=\frac{\left( \frac12 \right)_n^5\left( \frac12+k\right)_n^4\left( \frac15+\frac{2k}{5} \right)_n\left( \frac25+\frac{2k}{5} \right)_n\left( \frac35+\frac{2k}{5} \right)_n\left( \frac45+\frac{2k}{5} \right)_n\left( 1 + \frac{2k}{5}\right)_n}{(1)_n^5 (1+2k)_n \left( \frac12+\frac{k}{2}\right)_n^4\left( \frac12+\frac{k}{2}\right)_n^4} \frac{\left(\frac12 \right)_k^4}{(1)_k^4},
\]
and
\[
S(n,k)=\frac{640n^5(6n+4k+1)}{(n+2k+1)(5n+2k)}.
\]
From $F(n,k)$ we can get $R(n,k)$ such that $G(n,k)=U(n,k) R(n,k) z^n$, and using the technique explained in \cite{rama-divergent}, after replacing $n$ with $s$ and $k$ with $t$ we deduce that
\[
\frac{1}{2\pi i} \int_{-i \infty}^{i \infty} U(s,t) (1)_s \Gamma(-s) R(s,t) \left( \frac{3125}{1024} \right)^s \,ds = \frac{1280}{\pi^4} \quad \text{for all}\; t \in \mathbb{C},
\]
which gives the analytic continuation of
\[
\sum_{n=0}^{\infty} U(n,k) R(n,k) z^n, \quad |z|<1,
\]
at $z=-3125/1024$, that we can check with Maple for every complex number $k$ writing
\[
\texttt{sum(G(n,k),n=0..infinity,formal);}
\]
Specializing at $t=0$ gives
\begin{multline*}
\frac{1}{2\pi i} \int_{-i \infty}^{i \infty} \frac{\left( \frac12 \right)_s^5 \left( \frac15 \right)_s \left( \frac25 \right)_s\left( \frac35 \right)_s\left( \frac45 \right)_s}{(1)_s^8} \Gamma(-s) \left( \frac{3125}{1024} \right)^s \\ \times (5532s^4+5600s^3+2275s^2+425s+30) \, ds = \frac{1280}{\pi^4}.
\end{multline*}
In a completely similar way, from \cite[Example IX]{Pisco-seeds}, we obtain
\[
\frac{1}{2\pi i} \int_{-i \infty}^{i \infty} \frac{\left( \frac12 \right)_s \left( \frac15 \right)_s \left( \frac25 \right)_s\left( \frac35 \right)_s\left( \frac45 \right)_s}{(1)_s^4} \Gamma(-s) \left( \frac{3125}{256} \right)^s (483n^2+245n+30) \, ds = \frac{80}{\pi^2}.
\]

\section*{Ackowledgement}
I am very grateful to Wadim Zudilin for his comments and suggestions on this paper.

\end{document}